\documentclass[11pt,a4paper]{article}

\usepackage[margin=1in]{geometry}
\usepackage{amsmath, amssymb, amsthm}
\usepackage{graphicx}
\usepackage{subcaption}
\usepackage{float} 
\usepackage{booktabs}
\usepackage{array}
\usepackage{hyperref}
\usepackage{setspace}
\usepackage{authblk}
\usepackage{caption}
\usepackage{subcaption}
\usepackage{enumitem}
\usepackage{cite}
\usepackage{natbib}
\usepackage{longtable}

\theoremstyle{plain}

\theoremstyle{definition}

\theoremstyle{remark}

\title{\textbf{Synergies, Trade-offs, and Structural Pathways: A Directed Network Approach to SDG Prioritisation}}

\author{
  Gaurav Kottari\thanks{Corresponding Email: gk917@snu.edu.in},
  Niteesh Sahni\thanks{Email: niteesh.sahni@snu.edu.in}\\
  Department of Mathematics, Shiv Nadar Institution of Eminence,\\ Tehsil Dadri,  Gautam Buddha Nagar, 201314,\\ Uttar Pradesh, India 
}

\date{} 

\begin{document}

\maketitle

\begin{abstract}

To successfully implement the Sustainable Development Goals (SDGs), it is necessary to understand the process by which the achievement of one goal has a spillover effect in a development system. While existing research studies synergies and trade-offs among the SDGs, most empirical approaches operate at the goal level, treat interactions as undirected, or prioritise indicators without accounting for structural redundancy. In this paper, we propose a direction-sensitive and indicator-level network approach to detect high-impact and diversified entry points for policy intervention. By using statistically significant lagged correlations, we build a directed weighted network of SDG indicators and assign them into groups based on the balance of their positive and negative spillovers. Systemic effects are measured by weighted out-degree Opsahl centrality, and flow-based clustering is used to detect frequent paths of high positive spillovers. Applying the framework to the Indian context, it is found that the interlinkages in the SDGs are highly asymmetric and structured in specific structural subsystems. Although synergies slightly outweigh the total, trade-offs are still embedded in the sectors. Notably, the most influential indicators are focused on a single pathway of propagation, suggesting that influence-based prioritisation itself could result in redundant system-wide impacts. A cluster-based prioritisation approach leads to a more diversified set of interventions, triggering multiple structurally independent channels of beneficial spillovers. The proposed framework combines directionality, trade-off embedding, and structural propagation analysis in a single framework, providing a scalable solution for country-level SDG prioritisation under resource constraints.

\end{abstract}


\section{INTRODUCTION}
The 2030 Agenda for Sustainable Development provides a universal framework for advancing economic, social, and environmental progress simultaneously \citep{assembly2015resolution}. As a central component of this agenda, the Sustainable Development Goals (SDGs) articulate a collective vision for transformative progress through seventeen goals and their associated targets and indicators, adopted by all United Nations member states in 2015 \citep{le2015towards}. In this regard, the SDGs represent a shift from earlier international development initiatives, which often addressed different development sectors in isolation from one another. The SDG framework explicitly acknowledges that progress in one domain can influence outcomes in others \citep{nilsson2016policy, griggs2013sustainable}. This perspective reflects a broader shift within sustainability science toward systemic approaches that emphasise interactions and feedbacks across development domains \citep{stafford2017integration}. Consequently, sustainable outcomes emerge not only from progress within individual sectors but also from the interactions among them \citep{biermann2017global}.

While the SDGs emphasise integration and policy coherence, their implementation often follows existing institutional arrangements that remain largely sectoral in nature. Government agencies, development programmes, and policy instruments are typically organised around thematic domains such as health, energy, education, and environmental protection \citep{tosun2017governing, cejudo2017addressing}. As a result, many initiatives aimed at implementing the SDGs tend to focus on individual goals rather than addressing their interdependencies \citep{alcamo2020analysing}. Although this helps to make administration and management easier, it can lead to a lack of understanding of the overall implications of policy interventions. In this regard, it has been observed that policy interventions for the advancement of certain development goals can have unintended effects in other domains \citep{kroll2019sustainable, wong2019avoidance}. In this context, there has been increasing emphasis on the development of integrated policy approaches that take into account cross-domain interactions in sustainable development systems \citep{sachs2019six, breuer2019translating}.

In response to these challenges, research is increasingly focused on exploring and examining the interactions of SDGs. Empirical research carried out across various geographical locations suggests that interactions between different development goals may take the form of mutually reinforcing interactions and limiting interactions \citep{nilsson2016policy,pradhan2017systematic, tremblay2020sustainable}. For instance, it is possible that economic activities undertaken to improve health and education outcomes could also lead to positive economic and social outcomes, while economic activities could also impose negative pressure on the environment and impede progress toward ecological sustainability. The key point is that such interactions are likely to be context-dependent and could be influenced by different conditions and contexts of national development \citep{kottari2026network,bennich2020deciphering}. Therefore, it is observed that understanding how development interventions propagate through different and interconnected policy domains is a critical challenge for sustainability research and policy design \citep{scharlemann2020towards,alcamo2020analysing}.

Despite the recognition of these signed directional interlinkages, much of the existing literature employs an undirected and contemporaneous correlation-based approach to model SDG systems quantitatively. Centrality measures are commonly used to rank indicators while neglecting edge weights or considering all statistically significant relationships as equivalent \citep{swain2021modeling}. Moreover, trade-offs are usually not considered in the SDG network \citep{song2023unpacking}. Consequently, existing prioritisation strategies may miss three essential characteristics of SDG systems: directionality of influences, strength of spillovers, and simultaneous occurrence of both positive and negative outward influences. Without these factors, existing strategies may favor indicators that are important in an undirected correlation structure, but may not necessarily have a strong outward influence on SDG systems.

This limitation raises a central policy question: how can policymakers identify interventions that generate substantial positive spillovers while explicitly accounting for trade-offs and asymmetric influence patterns?

In order to bridge this gap, this study aims to create an indicator-level, direction-sensitive network approach for the analysis of interlinkages. Statistically significant lagged relationships are modelled as a directed weighted network, capturing time-ordered spillovers between indicators. The indicators are then categorised depending on the balance of their positive and negative out-linking effects. Opsahl weighted out-degree centrality \citep{opsahl2010node} is used for the calculation of systemic influence, where the magnitude of spillovers is taken into account.

In addition, the structure of the network is also investigated using a flow-based clustering technique that is particularly developed for directed weighted networks \citep{rosvall2008maps}. Unlike modularity-based community detection methods that identify the clusters of indicators mainly on the basis of the density of the edges, flow-based clustering identifies the groups of indicators on the basis of patterns of influence or the propagation of spillovers. From the viewpoint of policymakers, this can be considered as the propagation of spillovers from one indicator that can influence the other indicators that lie on the same influence pathways. The resulting clusters therefore represent groups of indicators that are connected through influence flows and may respond jointly to policy interventions.

Applying this framework to India reveals that the interlinkages between the SDGs are highly directional. Synergies and trade-offs coexist, and many indicators exert a strong outward influence on the broader system.  The application of flow-based clustering also indicates that groups of indicators are linked through specific influence pathways. Rather than focusing solely on key indicators in the entire system, it allows for the identification of representative indicators in each cluster, thus diversifying policy interventions over different subsystems in the SDG network.

This study contributes to the literature on interlinkages between SDG indicators in four ways. Firstly, it includes the concept of direction in indicator-level modelling. Secondly, it captures synergies and trade-offs in a unified framework. Thirdly, it uses weighted centrality to account for spillover magnitude. And finally, it proposes the concept of flow-based clustering, which groups indicators according to their connections in the form of influence propagation, facilitating a diversification-based prioritisation approach rather than relying on global centrality. By reframing SDG implementation as a problem of asymmetric and weighted influence flows, the paper advances both methodological precision and strategic guidance for integrated sustainable development planning.

The remainder of the paper is organised as follows. Section~\ref{4section:literature} reviews the existing literature on SDG interlinkages, distinguishing between qualitative and quantitative approaches and identifying key methodological limitations. Section~\ref{4section:methodology} presents the proposed network-based framework for modelling SDG indicator interactions, including the construction of the directed weighted network, the classification of indicators based on the balance of synergies and trade-offs, the measurement of systemic influence using weighted centrality, and the identification of structurally distinct propagation pathways through flow-based clustering. Section~\ref{4section:data} describes the SDG indicator dataset used in the empirical analysis and the preprocessing steps applied to construct a consistent time series suitable for network estimation.  Section~\ref{4section:results} reports the empirical results for India, including the structure of the directed SDG network, the classification of indicators, and the identification of high-impact synergistic indicators. Section~\ref{4section:discussion} interprets these findings in the context of national development policy and highlights their implications for SDG prioritisation strategies. Finally, Section~\ref{4section:conclusions} summarises the main contributions of the study and outlines directions for future research.

\section{LITERATURE REVIEW ON SDG INTERLINKAGES}\label{4section:literature}

\subsection{Qualitative approaches to SDG interlinkages}
Previous studies on interlinkages among SDGs have been based on qualitative methods, including expert opinion and judgement. Interaction scoring matrices and stakeholder-based evaluation frameworks have been used to classify interactions between goals and targets as synergistic or conflicting \citep{nilsson2016policy, griggs2017guide, allen2019prioritising, weitz2018towards, xiao2023synergies, fader2018toward}. These methods have served as an effective means to accomplish in proving that SDGs cannot be addressed in isolation, as well as creating awareness of inter-sectoral dependencies.

Despite their usefulness, qualitative methods have been shown to have limitations, including expert judgement, reproducibility, and scalability, especially in large sets of indicators \citep{pham2020interactions}. Although these methods help to create a better understanding of interlinkages, they offer little empirical support for the propagation of development over time \citep{allen2019prioritising}.

\subsection{Quantitative approaches to SDG interlinkages}

To address the limitations associated with qualitative scoring systems, a recent trend in the literature has been the application of quantitative methods for the analysis of SDG interlinkages. A significant share of the recent literature on SDG interlinkages uses statistical correlation methods for the analysis of synergies and trade-offs between SDG indicators. For example, \citet{de2020synergies} apply Spearman’s rank correlation method for the classification of interlinkages between SDG indicators in Spain using a set of pre-defined thresholds to distinguish between synergies and trade-offs. Similarly, \citet{kostetckaia2022sustainable} investigate the synergies and trade-offs between the SDGs across European Union countries and analyse how the share of synergies and trade-offs affects progress towards the 2030 Agenda using regression analysis.

In addition to pairwise correlation, some studies also employ dimension reduction methods to identify interactions at the goal level. \citet{hegre2020synergies} employ a dimension reduction method called Principal Component Analysis on the multi-dimensional structure of the SDG indicators and then examine the correlation between the components at the goal level across countries. Despite the reduction in arbitrariness in the choice of indicators and the improvement in the level of comparability, they still primarily reflect contemporaneous relationships rather than directional or structural propagation.

Extending correlation based analysis into network representations, \citet{miao2025priority} prioritise the SDGs at the goal level by integrating correlation coefficients, network analysis, and cluster analysis to determine implementation sequences in underdeveloped mountainous regions. Since policy interventions ultimately target specific indicators where measurable outcomes are achieved, several studies have developed indicator interaction networks to identify leverage points within the SDG system. In this context, \citet{swain2021modeling} have used interaction networks based on correlations, where centrality measures are used, including degree, eigenvector, and betweenness centrality, to identify influential indicators. However, the results of the study at aggregated regional levels show the presence of almost no strong negative correlations. On the contrary, at the country level, the SDG indicator network can have a high number of strong negative strength edges \citep{kottari2026network}. This contrast suggests that centrality-based prioritisation derived from predominantly positive interaction structures may not be directly applicable in national contexts where structural trade-offs are more significant.

Beyond correlation-derived networks, recent research has also employed semantic and systems-based approaches. For example, \citet{song2023unpacking} construct networks based on textual similarity among SDG targets using embedding techniques such as Word2Vec. While such approaches capture conceptual proximity, they do not explicitly incorporate trade-offs between targets and therefore provide only a partial representation of SDG interdependencies \citep{kottari2026network}.

Addressing these limitations, \citet{kottari2026network} proposed a network-based framework that changes the perspective from analysing the classification of individual interaction pairs to the analysis of the role of indicators within the interaction network. Instead of focusing on the nature of the relationships between the indicators as being either positive or negative, the proposed framework analyses the balance of the reinforcing and conflicting influences surrounding the individual indicators. The proposed analysis focuses on the relative importance of synergistic and trade-off relationships within the immediate neighborhood of the individual indicators. According to the proposed analysis, the indicators can be classified as synergy-dominated or trade-off-dominated within the interaction network. Indicators located within reinforcing structures have the potential to act as entry points for the policy intervention, whereas indicators located within trade-off structures require more coordinated intervention.

Despite this progress, systematic reviews continue to emphasise the need to improve dynamic modelling and the explicit inclusion of directional propagation mechanisms \citep{khot2026existing}. Although quantitative approaches have helped improve the empirical rigour and structural analysis, there is a need to further incorporate directional and dynamic modelling to effectively capture the complexity of SDG interdependencies.

In this context of the above limitations, this study takes a complementary approach to the above modelling paradigms. Rather than using pairwise associations, this study models SDG interlinkages at the indicator level as a directed weighted network, which explicitly includes reinforcing as well as conflicting relationships and can accommodate asymmetric influence structures in a unified framework.

To the best of our knowledge, this study is among the first attempts to quantitatively formalise the interdependencies between SDG indicators as a directed weighted network that captures the full range of interaction strengths. Importantly, our objective is not merely to identify a set of highly central indicators, as is common in much of the existing literature, but to examine the broader pattern of influences surrounding each indicator through its interactions with the remainder of the system. This perspective allows us to show that influential indicators are not homogeneous in their role within the SDG network; rather, they exhibit substantial diversity in their structural positions. Recognising this diversity provides a more refined understanding of potential policy entry points within the SDG framework.


\section{METHODOLOGY}\label{4section:methodology}
\begin{figure}[t]
    \centering
    \includegraphics[width=\textwidth]{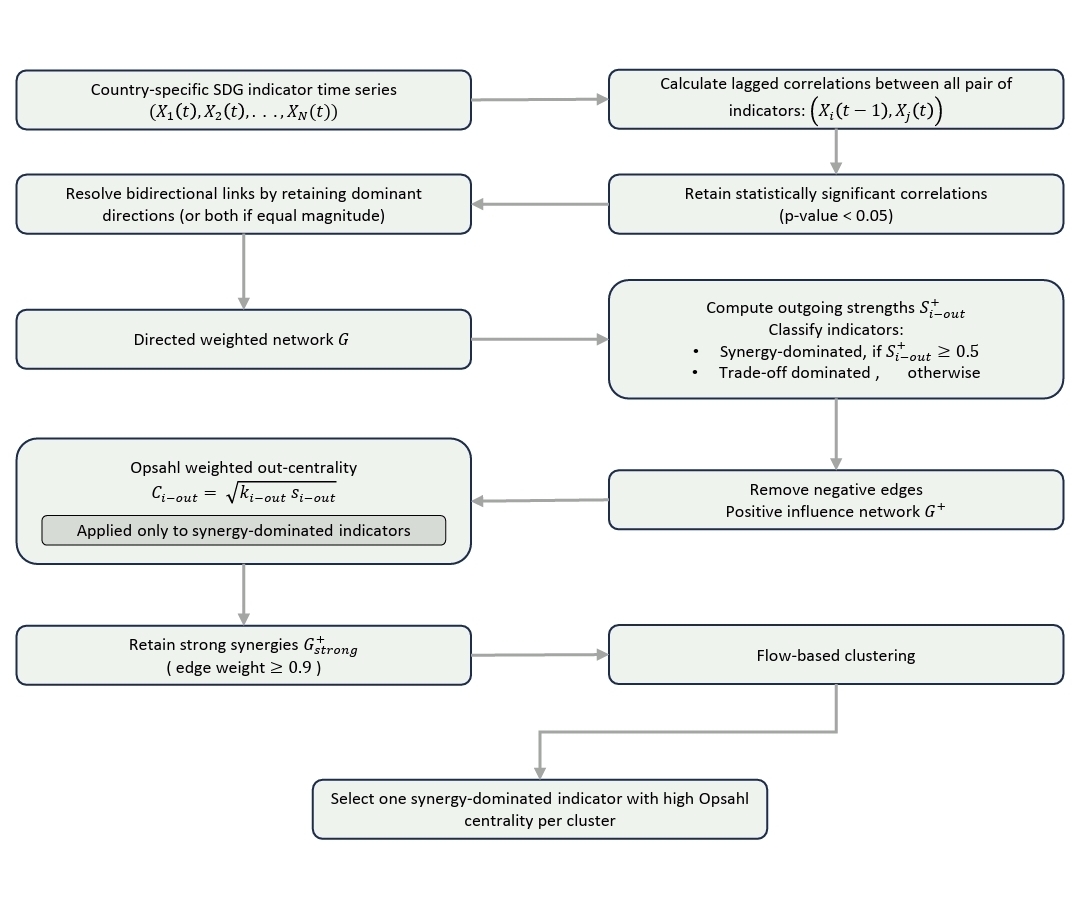}
    \caption{Overview of the methodological workflow for constructing directed SDG interlinkage networks, identifying synergy-dominated indicators, and prioritising high-impact indicators.}
    \label{4fig:methodology_flowchart}
\end{figure}

This section outlines a multi-step network-based framework to identify high-impact SDG indicators at the country level. Starting from country-specific SDG indicator time series, we construct a directed weighted network using lagged correlations, classify indicators according to the balance of synergies and trade-offs, quantify outward influence using weighted centrality, and finally identify structurally distinct pathways via flow-based clustering. Figure~\ref{4fig:methodology_flowchart} provides an overview of the methodological workflow.

\paragraph{Modelling SDG interlinkages using network theory.}

Sustainable development is inherently systemic, implying that gains made in one area have a tendency to spill over into other areas, either positively through reinforcing mechanisms or negatively through inhibiting mechanisms. Hence, when dealing with the SDG indicators for a particular country, it is possible that dealing with each indicator in isolation might miss out on the indirect effects, unintended consequences, and cumulative spillovers that are central to national development paths. In addition, it has been established that the nature and strength of the SDG interlinkages differ significantly between countries, depending on their economic structure, institutional arrangements, and development paths \citep{kottari2026network}. Network theory therefore provides a natural and rigorous framework for modelling country-specific SDG systems as interconnected structures rather than as independent targets.

In this paper, the SDG indicators are modelled as nodes, and the interlinkages between them are modelled as edges. Using this model, we can account for not only the presence of interdependencies between the SDG indicators but also their direction, weight, and sign, which are critical in the context of policy-relevant analysis. Hence, we model the system as a directed weighted network.


Formally, let
$$
G = (V, E, W)
$$
denote a directed weighted network, where $V = \{1,2,\dots,N\}$ represents SDG indicators, $E \subseteq V \times V$ the set of directed edges, and $W = \{w_{ij}\}$ the corresponding weights. A directed edge from indicator $i$ to indicator $j$ exists if $w_{ij} \neq 0,$ where $w_{ij}$ measures the impact of indicator \(i\) on indicator \(j\).

\paragraph{Estimating directional interlinkages and weights.}

To quantify the directionality of these interlinkages, we calculate the lagged correlations between the time series. This involves, for each pair of $(X, Y)$, estimating the correlation between $X_{t-1}$ and $Y_{t}$. This measures the delayed impact of changes in $X$ on changes in $Y$. Lagged correlations take into account time ordering by relating changes in a variable to future changes in another variable. The use of lagged correlations to infer directional influence follows approaches used in directed network analysis of time-dependent systems \citep{kullmann2002time, li2022undirected}. Only statistically significant correlations are retained as edges, using a hypothesis test with a significance threshold of $p\text{-value} < 0.05$, ensuring that the resulting network reflects robust interlinkages rather than noise-driven associations.  

In some cases, statistically significant bidirectional relationships may exist, i.e., $X_{t-1} \rightarrow Y_t$ and $Y_{t-1} \rightarrow X_t$. In order to maintain the parsimony of the network structure, we follow the dominance rule adopted in \citet{li2022undirected} for directed graphs. Here, the direction of the stronger relationship in terms of the absolute value of lagged correlation coefficients is preserved. This leads to the following rules: if $|w_{XY}| > |w_{YX}|$, only $X \rightarrow Y$ is preserved, if $|w_{YX}| > |w_{XY}|$, only $Y \rightarrow X$ is preserved, and if $|w_{XY}| = |w_{YX}|$, then both directions are kept.


\paragraph{Identifying synergy- and trade-off-dominated indicators.}

Once the directed weighted network is constructed, we classify indicators based on the balance of their contrasting outgoing effects. The proposed measures extend the positive and negative strength measures originally introduced for undirected weighted SDG networks \citep{kottari2026network}. Since the present framework models lagged relationships as a directed network, we focus specifically on the outgoing effects of each indicator. This is because of the policy-oriented nature of the analysis. Indeed, when a policy is implemented on a certain indicator, the interest is in understanding how the improvement in that indicator will propagate and affect other indicators in the system.

For each indicator \(i\), we define the positive and negative outgoing strengths as
\[
S^{+}_{i \text{-out}} = \frac{\sum_{j \neq i} \max(w_{ij},0)}
     {\sum_{j \neq i} |w_{ij}|},
\qquad
S^{-}_{i \text{-out}} = \frac{\sum_{j \neq i} |\min(w_{ij},0)|}
     {\sum_{j \neq i} |w_{ij}|}.
\]

An indicator is classified as \emph{synergy-dominated} if \(S^{+}_{i \text{-out}} \geq 0.5\), indicating that its net outward influence on the system is reinforcing. Conversely, indicators for which negative effects dominate are classified as \emph{trade-off-dominated}.

The major advantage of this classification is that it can measure the balance between reinforcing and inhibiting effects produced by an indicator. By combining positive and negative outgoing influences in one normalised framework, the proposed measures can determine whether an indicator has an overall outward effect that is dominated by synergies or trade-offs.

However, this classification has an important limitation: it identifies whether an indicator tends to generate synergies or trade-offs, but it does not quantify how broadly these effects spread across the system or how strongly they affect other indicators.  To understand this limitation, let us consider an example. Suppose that we are given that indicator \(i\) has three outgoing edges with weights -0.2, -0.2, and 0.8, whereas indicator \(j\) has four outgoing edges with weights 0.3, 0.4, 0.5, and -0.7. In this case, we can see that the positive outgoing strengths are \(S^{+}_{i \text{-out}} = 0.66\) and \(S^{+}_{j\text{-out}} = 0.63\). Therefore, both indicators are classified as synergy-dominated. However, their policy relevance differs substantially. Indicator $i$ exerts a strong positive influence on a single indicator, whereas indicator $j$ generates moderate positive influences on multiple indicators. Although both of these indicators are quite similar in their \(S^{+}_{i \text{-out}}\) values, they are quite different from each other.

This example shows that the classification based solely on the balance of indicators is not sufficient to find the high-impact indicators. From the point of view of the policy, the priority is to find indicators for which improvement has a significant impact on many other indicators. Nevertheless, the classification based on the synergy and trade-off is a necessary first step, as it allows us to exclude negative-weighted edges from subsequent analysis and focus explicitly on positive spillovers.

\paragraph{Measuring high-impact synergistic indicators.}

As discussed above, it is not sufficient to identify synergy-dominated indicators alone; we are primarily interested in indicators whose improvement is likely to generate the greatest downstream benefits. To quantify this impact, we employ the weighted out-degree centrality measure proposed by \citet{opsahl2010node}, which jointly captures the number of affected indicators and the cumulative strength of influence. For each indicator \(i\), Opsahl out-centrality is defined as follows,
\[
C^{\alpha}_{i\text{-out}} = k^{1-\alpha}_{i\text{-out}} \, s^{\alpha}_{i\text{-out}},
\]
where \(k_{i\text{-out}}\) denotes the number of outgoing connections from indicator \(i\), and \(s_{i\text{-out}}\) represents the total strength of its outgoing influence, computed as the sum of its outgoing edge weights. The parameter \(\alpha\) controls the relative importance of breadth versus intensity. When \(\alpha = 0\), only the number of outgoing links matters; when \(\alpha = 1\), only their total weight matters. Since effective SDG interventions require both wide reach and substantial impact, we set \(\alpha = 0.5\), assigning equal importance to both dimensions. Under this choice, the Opsahl out-centrality reduces to
\[
C_{i\text{-out}} = \sqrt{k_{i\text{-out}} \, s_{i\text{-out}}}.
\]

This measure is well-defined only when \(s_{i\text{-out}} \geq 0\). Accordingly, after classifying indicators using $S^+_{i\text{-out}}$, we remove all negative-weighted edges and define the positive influence network
\[
G^{+} = (V, E^{+}, W^{+}),
\quad
E^{+} = \{(i,j) \in E : w_{ij} > 0\}.
\]

 Opsahl out-centrality is then computed only for synergy-dominated indicators in $G^+$, as our objective is not merely to identify influential indicators, but to prioritise among those whose outward influence is predominantly influential and thus aligns with the objectives of sustainable development.

\paragraph{Strengthening policy relevance through strong-synergy structures.}

While the ranking of indicators based on the Opsahl out-centrality score alone provides a measure of the strength and reach of the influence of an indicator on others, it does not take into consideration the propagation of the effect of the intervention from one indicator to another. In reality, the improvement of an indicator initiates a chain of interactions with other indicators, which are closely interconnected and form cohesive impact pathways. For instance, an investment in increasing the share of renewable energy propagates through improvements in electricity access, which subsequently influence years of schooling and ultimately life expectancy \citep{collste2017policy},
\[
\text{Renewable energy share}
\;\longrightarrow\;
\text{Electricity access}
\;\longrightarrow\;
\text{Years of schooling}
\;\longrightarrow\;
\text{Life expectancy}.
\]
Similarly, poverty reduction interventions may propagate through a largely overlapping pathway,
\[
\text{Poverty reduction}
\;\longrightarrow\;
\text{Years of schooling}
\;\longrightarrow\;
\text{Life expectancy}.
\]

The examples above illustrate that different policy interventions can have similar effects on intermediate and downstream indicators. This means that, in effect, it is possible that prioritising all highly influential indicators might result in redundant policy interventions as their effects cascade down similar pathways. It is also important to note that improving all indicators is desirable and consistent with the holistic nature of sustainable development. However, there are many financial, institutional, and administrative constraints that make it difficult for countries to address all indicators simultaneously. In that respect, policy prioritisation is no longer about ranking indicators based on their influence but rather about finding a subset of indicators that have the most comprehensive coverage in terms of diversity. This is where the need arises to consider not only the influence of each indicator but also how influence is structured and cascaded in the network.

\begin{figure}[h]
    \centering
    \includegraphics[width=1\textwidth]{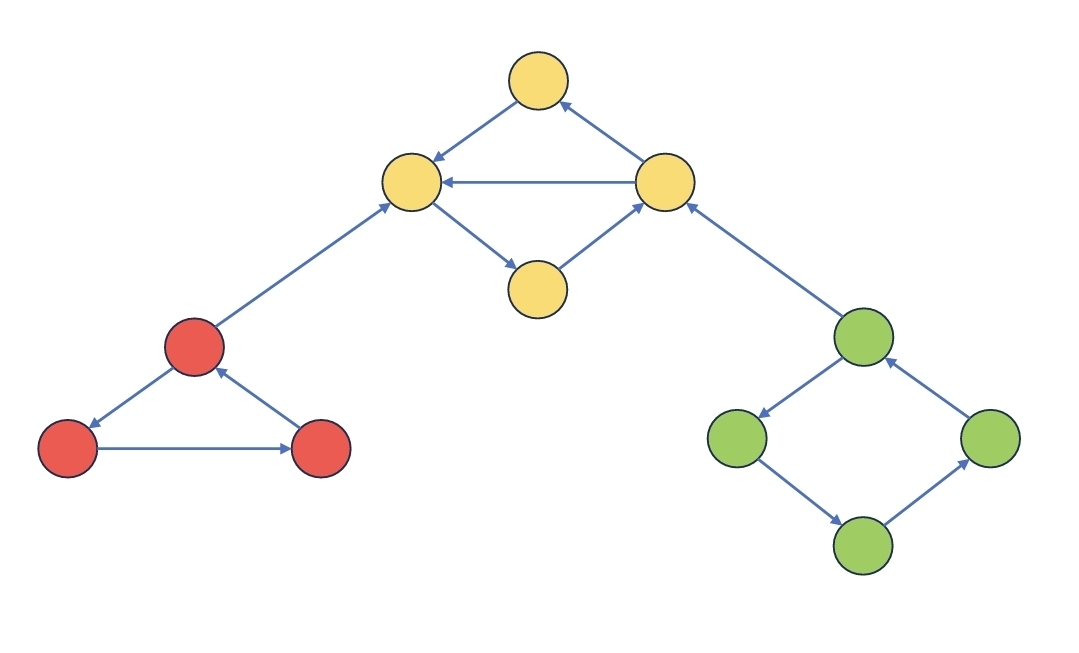}
    \caption{Schematic illustration of flow-based clustering. Nodes of the same colour indicate clusters identified based on recurrent influence flows, highlighting coherent pathways of propagation rather than simple connection density.}
    \label{4Flowbasedclusteringillustration}
\end{figure}

For the purpose of capturing such a structure, we use a method called \emph{flow-based clustering} \citep{rosvall2008maps}, where indicators are grouped based on the patterns of influence propagation within the network. The clustering is computed using the Infomap implementation from the MapEquation software package \citep{mapequation2026software}. This method identifies groups of indicators among which influence tends to propagate and persist more frequently than with the rest of the system. For instance, as depicted in Figure~\ref{4Flowbasedclusteringillustration}, the indicators that share the same colour are grouped as clusters based on the fact that there is a coherent flow within the group. Within each coloured group, influence propagates along directed sequences of links and tends to remain within the group for successive transitions of the network flow, while movements across differently coloured groups occur less frequently. The clusters are therefore subsystems within the network through which the influence of the policies tends to propagate. The fact that the method takes into account the direction and strength of the interlinkages makes it more appropriate for identifying structurally distinct subsystems within the network.

We apply this clustering approach to a reduced network \(G^{+}_{\text{strong}}\), constructed from \(G^{+}\) by retaining only positive interlinkages with edge weights \(w_{ij} \geq 0.9\). Although all links in \(G^{+}\) are statistically significant, the clusters formed over \(G^{+}\) might not capture the major channels of influence due to the presence of weak as well as strong spillovers. Clustering based on strong synergies would therefore highlight the major channels through which strong positive spillovers are transmitted in the network. Further, form each cluster, we then choose the synergy-dominated indicator with the highest Opsahl out-centrality. Instead of choosing the top-ranked indicators across the entire system, we aim to target unique structural elements of the system, thereby ensuring the maximum impact while minimising redundancy. 

\section{Data and Preprocessing}\label{4section:data}

The SDG indicator dataset created by the SDG Transformation Center and published through the Sustainable Development Report (SDR) database is used in this study's empirical analysis.\footnote{Data accessible at \url{https://dashboards.sdgindex.org/downloads/}.} Using information from numerous international statistical sources, the database offers yearly country-level time series for a wide range of indicators related to the SDGs. The SDR dataset includes both raw indicator values and transformed “backdated” indicators used in the construction of the SDG Index.

Instead of using raw data in this study, we used backdated dataset. Each indicator in the backdated series is rescaled to a common range between 0 and 100 using a min-max normalisation framework \citep{sachs2025sdr}. In this case, reaching the SDG target is represented by 100, while the lowest performance is represented by 0. This modification ensures that improved performance across all indicators is demonstrated by higher values. Because the direction of improvement varies among raw SDG indicators, standardisation is crucial. Specifically, for some indicators, higher values represent progress (e.g., access to services), whereas for others lower values represent progress (e.g., poverty or mortality rates). By standardising the direction of progress across indicators, the transformed dataset enables correlations to be consistently interpreted as representing either synergies or trade-offs.

The data set ranges from 2000-2024 for most indicators. However, the construction of the network is based on lagged correlations using a one-year lag. Hence, indicators that have missing data during this period are removed. Also, indicators that have constant values during the period from 2000-2023 or from 2001-2024 are removed because they have no variation, and hence no information for the estimation of the interdependencies. After applying the filters, the remaining indicators will have enough information for the construction of the directed weighted SDG interlinkage network.

\section{RESULTS}\label{4section:results}
\subsection{Directed SDG interlinkage network}
Using the methodology that has been proposed in this paper on the data from India, the resulting directed and weighted SDG indicator network comprises 79 indicators and 2,643 statistically significant lagged interlinkages. The list of the 79 SDG indicators is presented in Table~\ref{4indicatorslist} (see Appendix).

\begin{figure}[t]
    \centering
    \includegraphics[width=0.85\textwidth]{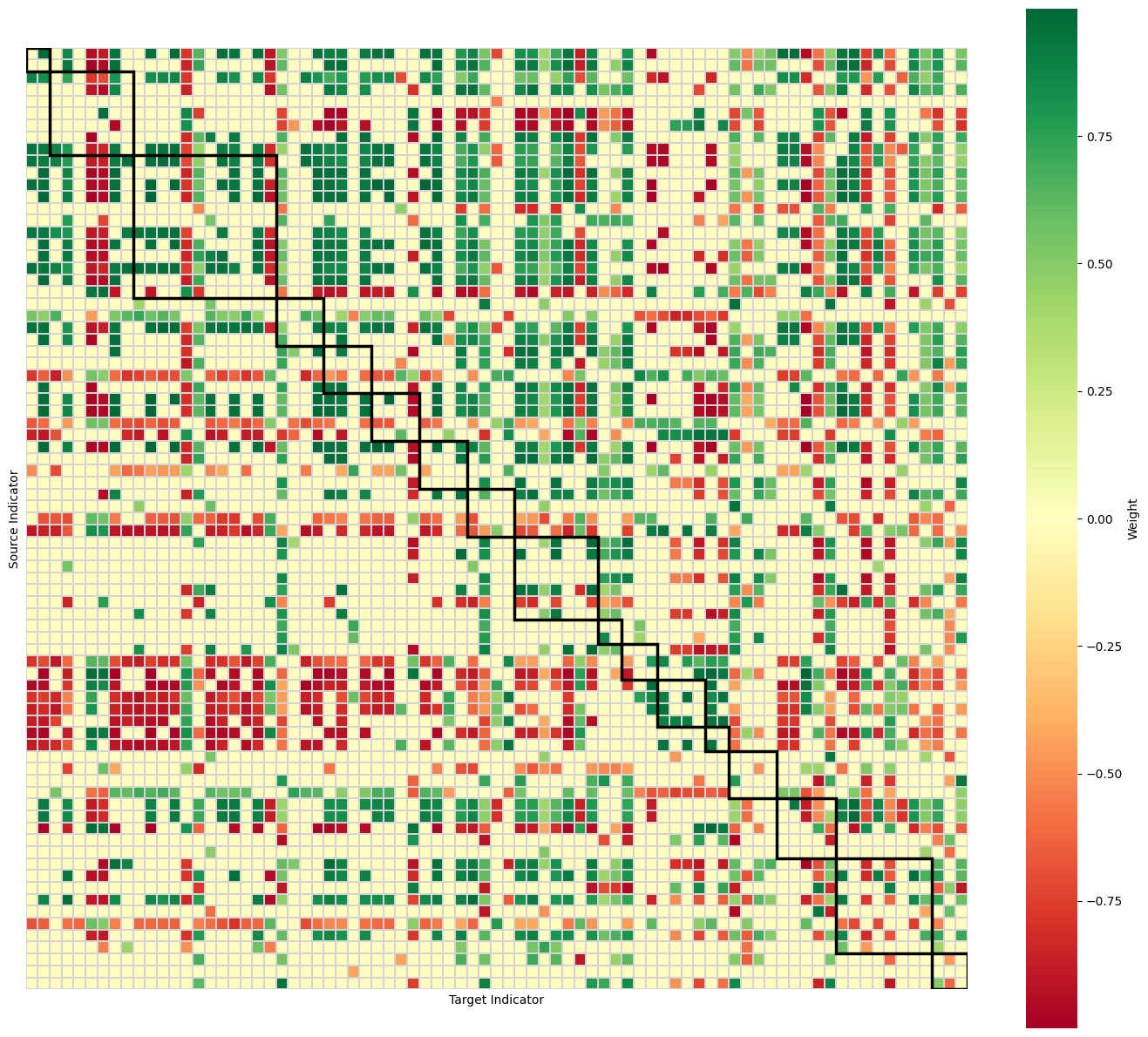}
    \caption{Heat map of statistically significant dominant lagged correlations among SDG indicators. Rows denote source indicators and columns denote target indicators. Indicators are ordered by SDG, with diagonal blocks corresponding to within-goal interlinkages (SDG 1 at top-left to SDG 17 at bottom-right). Off-diagonal blocks capture cross-goal indicator influences.}
    \label{4fig:lagged_heatmap}
\end{figure}

The SDG network for India shows that positive spillovers prevail in the system. It has been observed that 1,491 edges (56.4\%) of the Indian SDG network have positive weights, while the remaining 1,152 edges (43.6\%) of the network depict negative interlinkages. Although positive edges prevail in the SDG system, the existence of significant negative edges shows that trade-offs have remained an integral part of the Indian SDG system. A summary of the lagged correlations underlying this network is provided in Figure~\ref{4fig:lagged_heatmap}. Here, diagonal blocks represent the interdependencies among indicators within the same goal, whereas numerous entries appear in the off-diagonal blocks, representing interdependencies among indicators belonging to different goals. However, it should be noted that the impacts of indicators on the same or different goals can be both positive and negative.

In general, these patterns show that the SDG interlinkages in India are neither confined to individual goals nor uniformly reinforced. This provides empirical evidence on why the indicator-level classification based on the balance of synergies and trade-offs in outward influence is essential.

\subsection{Classification of indicators}
When classifying indicators through synergy-trade-off classification using the measure $S^{+}_{i \text{-out}}$, significant heterogeneity emerges with regard to the systemic role performed by each indicator within the overall structure. Of the total of 79 indicators considered in the study, 52 fall into the category of indicators that are synergy-dominated, while 27 fall into the category that are trade-off-dominated indicators. This indicates that for most indicators, positive downstream effects outweigh negative ones. At the same time, the non-negligible share of trade-off-dominated indicators highlights that inhibiting effects remain structurally embedded within the Indian SDG system.

\begin{figure}[t]
    \centering
    \includegraphics[width=0.65\textwidth]{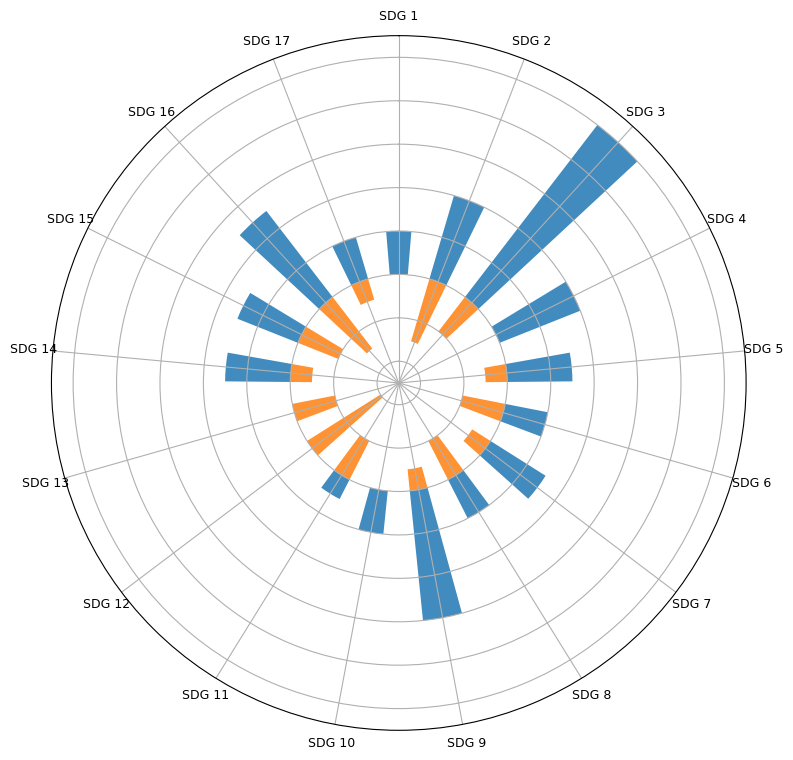}
    \caption{Polar representation of synergy- and trade-off-dominated indicators across the SDGs. For each SDG, outward (inward) radial bars represent the number of synergy- (trade-off-) dominated indicators. The radial scale is discretised such that each unit of bar length corresponds to two indicators; thus, bar length directly encodes indicator counts.}
    \label{4fig:synergydominateddistribution}
\end{figure}

The distribution of the synergy-dominated indicators with respect to the goals is not uniform, and it is evident from Figure~\ref{4fig:synergydominateddistribution}. The highest number of synergy-dominated indicators has been recorded in goal 3 (Good Health and Well-being). The second highest is in goal 9 (Industry, Innovation, and Infrastructure). This indicates that improving in the domain of health and innovation-related factors will lead to positive effects in almost all development domains and hence will be structurally significant entry points in the system. In contrast, all indicators studied related to goal 12 (Responsible Consumption and Production) and goal 13 (Climate Action) were identified as trade-off-dominated. This reflects that their progress is associated with a net negative downstream influence on other indicators. Thus, environmental and production-related interventions should be accompanied by supporting actions to prevent negative effects from spreading to other development outcomes.

A particularly illustrative case is the indicator: Access to improved piped water source. Although intuitively this may be considered a development priority, the proposed framework classified this indicator as trade-off-dominated $(S^{+}_{i \text{-out}}=0.276)$. In fact, it features six outgoing interlinkages with health-related indicators under SDG 3, of which five are strongly negative. Four of these negative effects have magnitudes close to $-0.9$, while the remaining one is close to $-0.6$. This implies that an increase in access to piped water is followed by deteriorations in several health outcomes. This pattern suggests that expansion of the infrastructure alone, when not accompanied by adequate control over water quality or sanitation, may generate adverse unintended consequences. 

In contrast, several indicators are classified as synergy-dominated. These include fish caught by trawling or dredging, corruption perceptions Index, renewable energy share in total final energy consumption, literacy rate, and prevalence of undernourishment. This indicates that synergy-dominated indicators are present across multiple domains, including environmental management, governance, energy, education, and food security.

\subsection{Identification of high-impact synergistic indicators}
Having identified synergy-dominated indicators and constructed the positive influence network $G^+$, the next step is to determine which indicators are influenced to the greatest extent. According to the methodology, this is achieved by employing the Opsahl weighted out-centrality measure, which jointly captures the breadth of influence (number of affected indicators) and the cumulative strength of that influence.

The ten highest-ranked synergy-dominated indicators according to Opsahl out-centrality are listed in Table~\ref{4tab:opsahl_top10}. These include surviving infants who received two WHO-recommended vaccines, lower secondary completion rate, poverty headcount ratio at \$2.15 per day, under-five mortality rate, sustainable nitrogen management index, maternal mortality ratio, population with access to electricity, population using at least basic drinking water services, adolescent fertility rate, and the corruption perceptions index. They cover all important aspects of development, such as health, education, poverty, governance, infrastructure development, and environmental management.

\begin{table}[t]
\centering
\caption{Top ten synergy-dominated SDG indicators ranked by Opsahl out-centrality.}
\label{4tab:opsahl_top10}
\begin{tabular}{clrr}
\hline
SDG & Indicator & Opsahl out-centrality \\
\hline
SDG 3 & Surviving infants who received two WHO-recommended vaccine 
& 38.57 \\
SDG 4 & Lower secondary completion rate  
& 37.92 \\
SDG 1 & Poverty headcount ratio at \$2.15/day 
& 36.84 \\
SDG 3 & Mortality rate, under-5 (per 1,000 live births) 
& 35.64  \\
SDG 2 & Sustainable Nitrogen Management Index  
& 34.98 \\
SDG 3 & Maternal mortality ratio (per 100,000 live births) 
& 34.48 \\
SDG 7 & Population with access to electricity 
& 32.71 \\
SDG 6 & Population using at least basic drinking water services 
& 32.56 \\
SDG 3 & Adolescent fertility rate (births per 1,000 females aged 15--19) 
& 32.43 \\
SDG 16 & Corruption Perceptions Index 
& 32.29 \\
\hline
\end{tabular}
\end{table}

To understand the mechanism through which positive development effects flow through the SDG system, it is essential to begin by analysing the existing pattern of strong synergistic relations through flow-based clustering on the strong positive influence network \(G^{+}_{\text{strong}}\), considering only interlinks that have at least 0.9 in terms of strength.

\begin{figure}[h!]
\centering

\begin{subfigure}{0.48\textwidth}
    \centering
    \includegraphics[width=\linewidth]{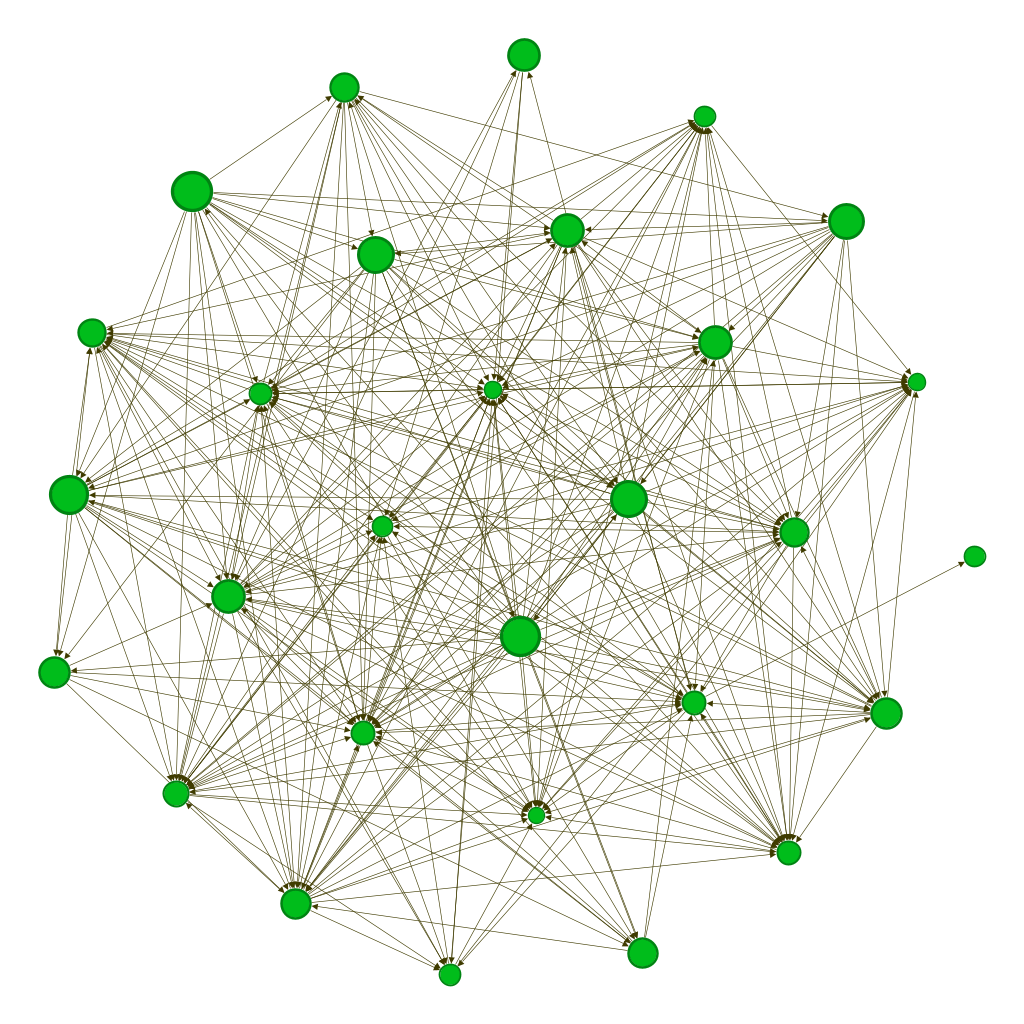}
    \caption{Cluster 1}
\end{subfigure}
\hfill
\begin{subfigure}{0.48\textwidth}
    \centering
    \includegraphics[width=\linewidth]{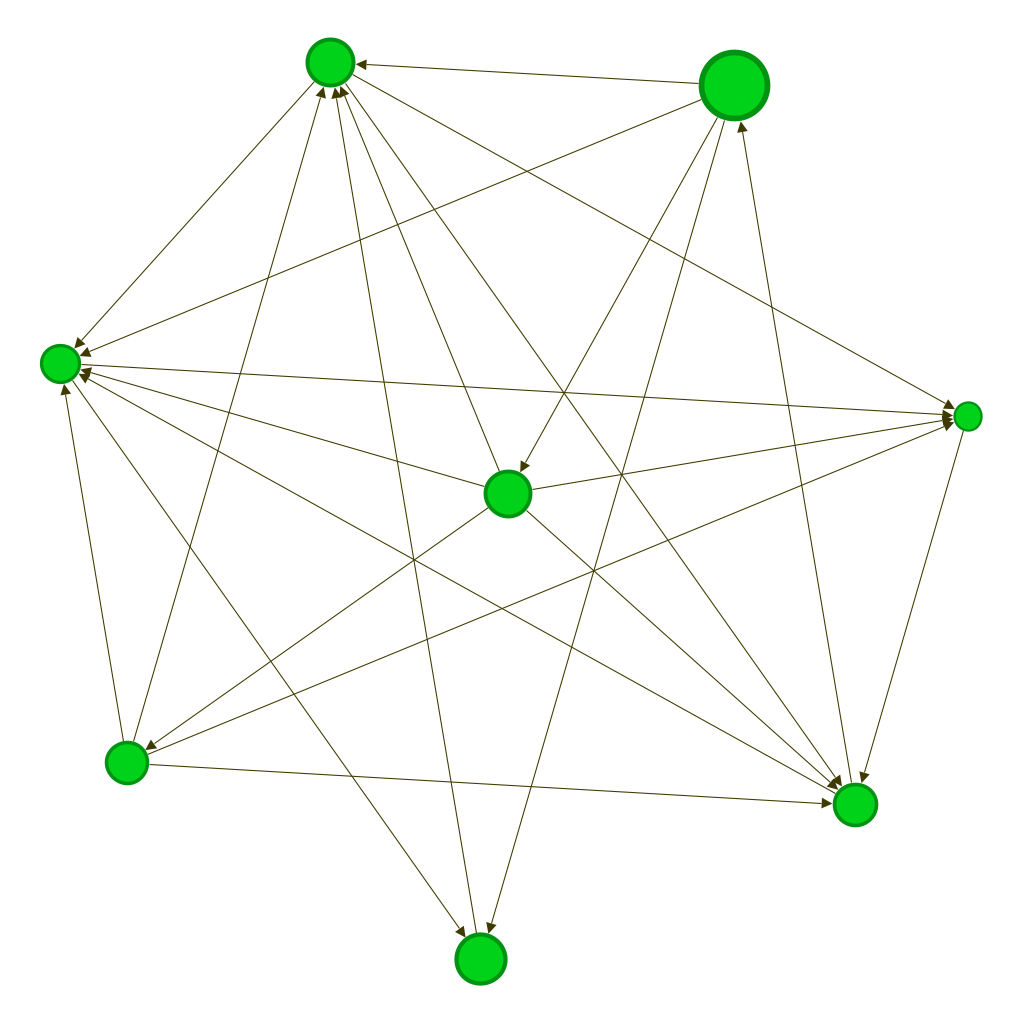}
    \caption{Cluster 2}
\end{subfigure}

\vspace{0.4cm}

\begin{subfigure}{0.48\textwidth}
    \centering
    \includegraphics[width=\linewidth]{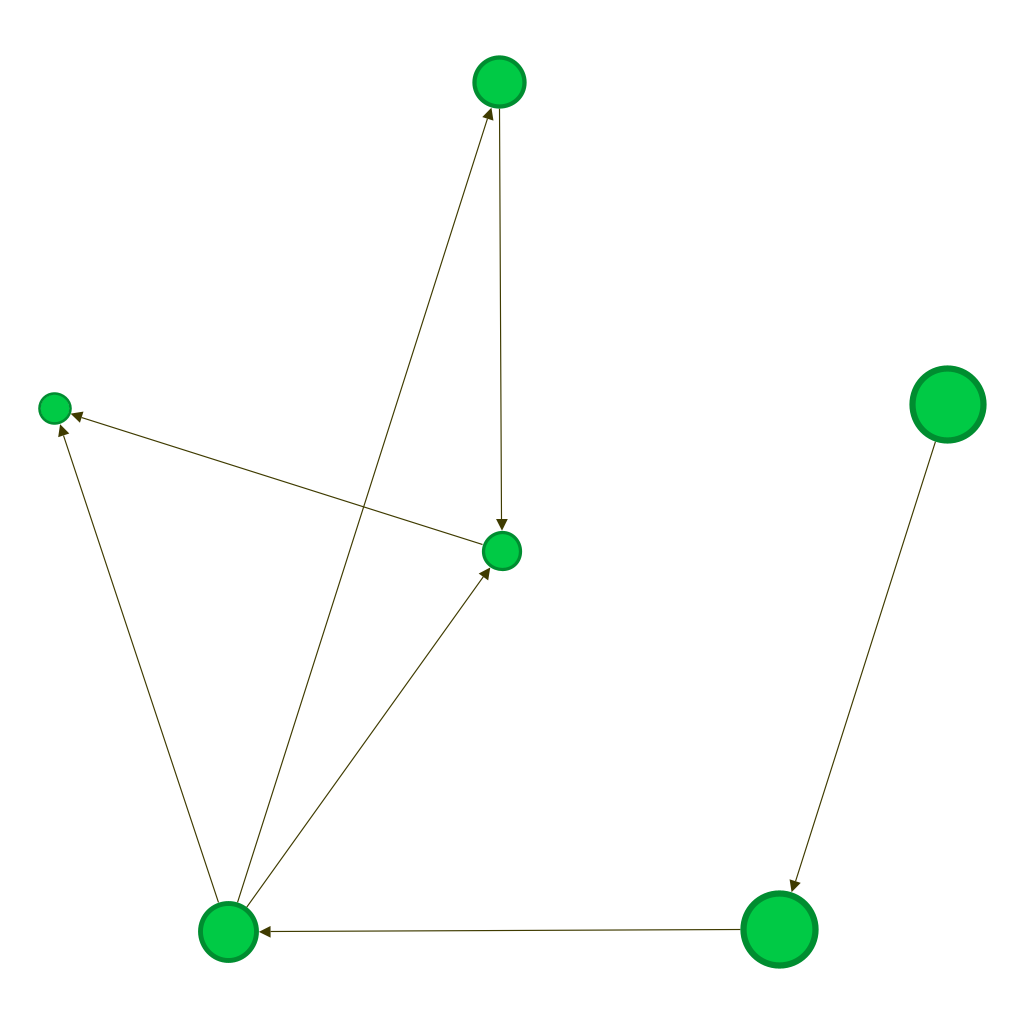}
    \caption{Cluster 3}
\end{subfigure}
\hfill
\begin{subfigure}{0.48\textwidth}
    \centering
    \includegraphics[width=\linewidth]{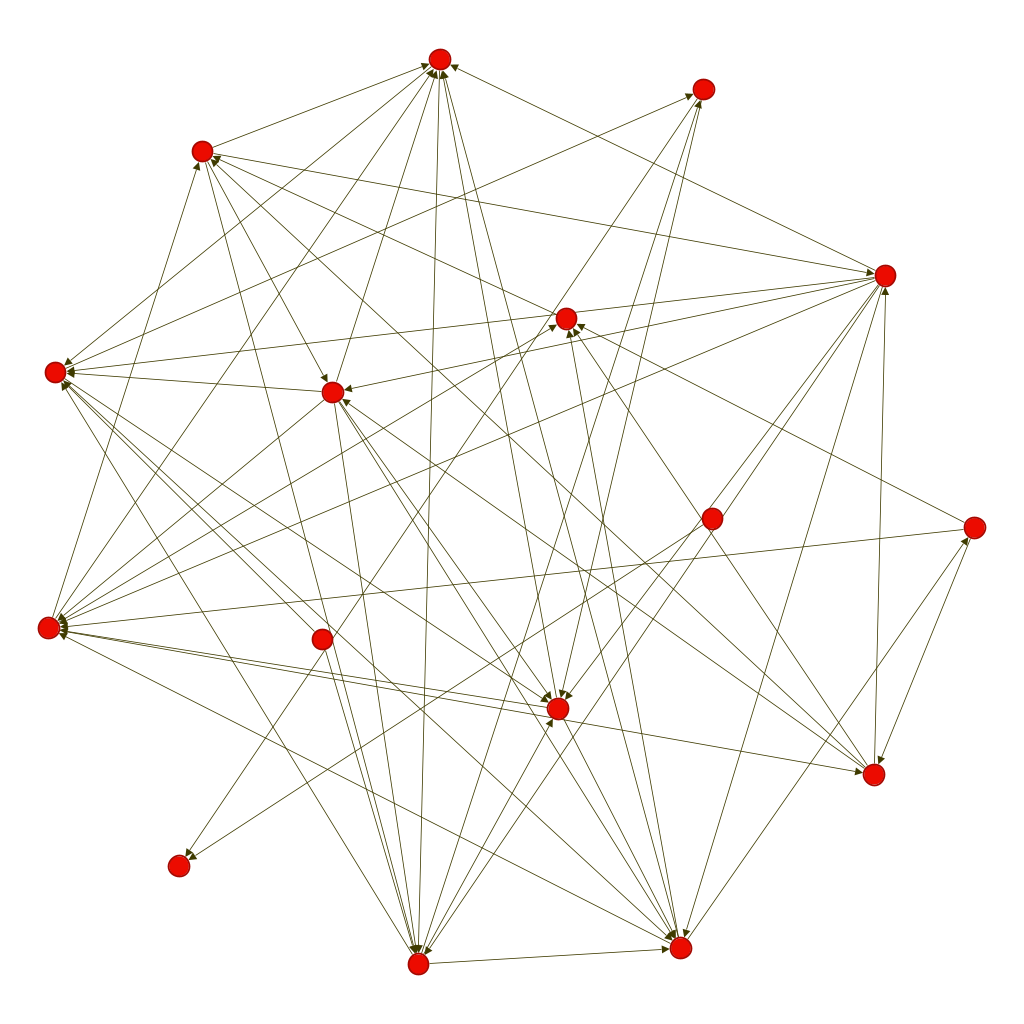}
    \caption{Cluster 4}
\end{subfigure}

\caption{Flow-based multi-indicator clusters of the strong positive influence network. Green nodes represent synergy-dominated indicators and red nodes represent trade-off-dominated indicators. Among the green nodes, node size is proportional to Opsahl out-centrality and is scaled separately within each cluster; sizes are not comparable across clusters.}
\label{4fig:fourclusters}

\end{figure}

Flow-based clustering resulted in four multi-indicator clusters (Figure~\ref{4fig:fourclusters}) and twenty single-indicator clusters. The multi-indicator clusters represent structurally cohesive subsystems in which strong synergistic effects circulate among several indicators, whereas the single-indicator clusters correspond to indicators whose strong positive effects do not form closed propagation structures with others.

Among the various multi-indicator clusters, Cluster 2 represents a structurally integrated synergy-dominated subsystem with indicators such as gender parity in education, share of renewable energy, internet use, mobile broadband subscriptions, performance in university ranking, publication output, patent applications, and prevalence of urban slum. All these indicators have strong interconnections among themselves in terms of directed lagged links. The presence of directed cycles also provides evidence of reinforcing feedback in the subsystem, which represents the systematic propagation of changes through paths and back to affect earlier nodes. Cluster 2, therefore, captures a development configuration that is highly integrated with strong lagged synergy interdependencies in education, energy transition, digital infrastructure development, innovation performance, and urban development.

In contrast to Cluster 2, Cluster 1 forms a socio-developmental subsystem in which strong positive spillovers circulate across poverty, health, education, and governance indicators. Cluster 3 connects social, institutional, and environmental indicators, forming a smaller but coherent pathway of influence. Cluster 4 captures an environmental–resource subsystem in which biodiversity, emissions, resource use, and pollution-related indicators interact through strong synergistic feedbacks. Together, these clusters show that strong positive spillovers are organised within specific groups of related indicators rather than influencing the system in a uniform manner.

A notable finding is that all ten indicators in Table~\ref{4tab:opsahl_top10} belong to Cluster~1. The strong concentration of these top-ranked indicators within a particular cluster reveals that their spillovers operate through largely similar channels. Consequently, prioritising multiple top-ranked indicators from this cluster would generate partially redundant system-wide effects, as improvements tend to reinforce similar intermediate and downstream outcomes.

Cluster-based prioritisation strategies presented in Section~\ref{4section:methodology} produce a structurally diversified set of high-impact indicators with predominantly positive spillover effects. This is achieved by selecting the synergy-dominated indicator with the highest Opsahl out-centrality from each cluster. Consequently, Cluster~1 yields ``Surviving infants who received two WHO-recommended vaccines", Cluster~3 yields the ``Statistical Performance Index”, and Cluster~2 yields ``Articles published in academic journals”. Cluster~4 does not contribute an indicator, as none of its members are classified as synergy-dominated. In addition, all isolated synergy-dominated indicators are included in the prioritised list. Here, ``isolated synergy-dominated indicators'' refer to indicators that form single-indicator clusters in the strong positive influence network $G^{+}_{\text{strong}}$ but are nevertheless classified as synergy-dominated according to $S^{+}_{i\text{-out}} \ge 0.5$. The resulting prioritised set is reported in Table~\ref{4tab:cluster_prioritised}, and it differs substantially from the ranking based solely on Opsahl out-centrality shown in Table~\ref{4tab:opsahl_top10}.

\begin{table}[t]
\centering
\caption{Cluster-based prioritised synergy-dominated SDG indicators.}
\label{4tab:cluster_prioritised}
\begin{tabular}{cll}
\hline
Category & SDG & Indicator \\
\hline
Cluster 1 & SDG 3 & Surviving infants who received two WHO-recommended vaccines (\%) \\
Cluster 2 & SDG 9 & Articles published in academic journals (per 1,000 population) \\
Cluster 3 & SDG 17 & Statistical Performance Index (0--100) \\
Isolated & SDG 2 & Prevalence of undernourishment (\%) \\
Isolated & SDG 3 & Traffic deaths (per 100,000 population) \\
Isolated & SDG 4 & Net primary enrollment rate (\%) \\
Isolated & SDG 9 & Logistics Performance Index: Infrastructure score (1--5) \\
Isolated & SDG 10 & Gini coefficient \\
Isolated & SDG 10 & Palma ratio \\
Isolated & SDG 14 & Fish caught that are then discarded (\%) \\
Isolated & SDG 16 & Expropriations are lawful and adequately compensated (0--1) \\
Isolated & SDG 17 & Government spending on health and education (\% of GDP) \\
\hline
\end{tabular}
\end{table}

\section{DISCUSSION}\label{4section:discussion}

The results of the empirical analysis reveal that the system of SDGs in India is a complex web of directional spillovers between its indicators. While positive spillovers are predominant, a significant number of trade-offs are embedded in the system, implying that development progress comes with both positive and negative spillovers. The results further show that influential indicators are concentrated within specific structural subsystems, highlighting the importance of considering propagation pathways rather than relying solely on indicator rankings. In the following discussion, we interpret these findings in terms of policy coordination, management of trade-offs, and the identification of strategic entry points to accelerate SDG progress.

\paragraph{Directional spillovers and the limits of goal-based policy design.}

The directed lagged network highlights that the process of development in India is asymmetric and does not depend only on individual goals. There are a large number of statistically significant interlinks between different goal indicators, and these interlinks are directional in nature. This result further strengthens the point that the design of policies purely based on the 17-goal framework might obscure the true channels of development effects.

The strong presence of off-diagonal blocks in the lagged correlation heat map (Figure~\ref{4fig:lagged_heatmap}) suggests that better results in one area tend to precede changes in other areas. This implies that the sequencing and coordination of policies are important. Policies designed in isolation might underestimate the delayed effects, both positive and negative, that spread through the system as a whole. Ministries and departments working in their narrowly defined sectors might inadvertently create spillovers that reveal themselves only after a delay. The network-based evidence therefore supports more integrated planning frameworks that explicitly account for inter-ministerial dependencies.

\paragraph{Managing trade-offs in a predominantly synergistic system.}

Although the number of positive spillovers is slightly higher than that of negative ones in the Indian SDG system, the magnitude and distribution of trade-offs are far from negligible. The fact that 27 out of 79 indicators are trade-off-dominated indicators clearly shows that development progress can have inhibiting factors if not properly managed.

The case of piped water access provides a particularly illustrative example. Although considered a fundamental development area, the indicator has strong negative effects on several health-related indicators. This finding does not mean that the expansion of piped water access is not a good development target; it only means that the lack of proper protection mechanisms, such as water quality and sanitation, can have unforeseen consequences. Although previous studies~\citep{kottari2026network} found that piped water access is trade-off-dominated in the Indian SDG system, the current directed network analysis framework provides further evidence on the direction of influence rather than mere co-movements. 

Similarly, the existence of trade-off-dominated indicators in Goals 12 and 13 indicates that environmental and climate change transitions might create pressures in other parts of the system. This further implies that better performances in environmental indicators are succeeded by negative changes in socio-economic outcomes. This finding is supported by empirical studies showing that environmental and resource policies create short-term adjustment pressures on economic growth, employment, or productivity before long-term benefits are achieved. For instance, cross-country empirical findings show the existence of adverse effects of stricter environmental policies on productivity in the short term~\citep{albrizio2014empirical}. In addition, empirical findings in China show that environmental policies might affect employment outcomes through the influence of industry restructuring and compliance costs~\citep{cao2017effect}. Integrated models of water resource management show that environmental constraints might result in economic and employment losses despite improvements in environmental outcomes~\citep{zhao2021quantifying}. Therefore, environmental ambition should be embedded within coordinated transition strategies that mitigate adverse effects. 

In general, these findings indicate that progress in the infrastructure or production-related domains should be supported by complementary policies to avoid unintended negative spillovers.

\paragraph{Health, education, and governance as systemic entry points.}
The ranking based on Opsahl out-centrality reveals a set of indicators, mostly in the domain of health, education, poverty reduction, and governance, which have strong positive downstream effects. The dominance of health indicators in the top-ranked indicators highlights the positive effect of improving maternal health, child mortality, and vaccination, which goes beyond the domain and provides a foundation for other aspects of development.

This is consistent with the larger body of evidence that investing in people in terms of their health, skills, and overall well-being is critical to sustainable development. In fact, studies have consistently demonstrated that healthier people result in higher productivity and therefore higher income and overall economic performance ~\citep{bloom2004effect}. Moreover, studies on the interlinkages of SDGs reveal that the health and education goals have significant positive spillovers to other key areas of development ~\citep{pradhan2017systematic}. Therefore, in improving child and maternal health, it is not just about improving their well-being but also their educational attainment, productivity, and earning capacity, which in turn contributes to poverty reduction and stability.

However, what is important to note from the results of the clustering is that all ten of the top-ranked indicators according to Opsahl out-centrality are part of Cluster 1. This means that the effect of these indicators on the system is primarily channelled through overlapping socio-developmental tracks. While each of these indicators is highly significant on its own, prioritising these indicators at once may yield diminishing marginal returns in terms of system coverage because of their overlap on intermediate nodes.

\paragraph{Structural clustering and the avoidance of policy redundancy.}
Flow-based clustering reveals the existence of four different subsystems in which strong synergies circulate. In contrast to the modularity-based clustering methods applied in previous studies in SDG networks \citep{swain2021modeling,song2023unpacking}, the current method distinguishes clusters based on recurrent directional flows of influence. Modularity-based methods are capable of capturing high levels of connectivity but are not necessarily indicative of the influence of effects on the system. In contrast to modularity-based methods, flow-based clustering is able to capture the direction and strength of the interactions and identify the pathways through which positive effects are propagated, such as feedback cycles. From a policy point of view, this is an important difference, as it is not sufficient to identify groups of densely connected indicators, as this does not necessarily imply that interventions in these groups would be able to activate mechanisms of sustained propagation, while flow-based clustering directly identifies structurally coherent channels of improvement.

The formed clusters thus reflect meaningful influence pathways. Cluster 1 captures a socio-developmental pathway relating to poverty, health, education, and governance. Cluster 2 encompasses digital connectivity, knowledge creation, renewable energy transition, and urban development. Cluster 3 represents a smaller path related to institutional and social indicators. Cluster 4 represents an environmental-resource subsystem, which is defined by strongly interacting ecological processes. The concentration of highly ranked indicators in Cluster 1 thus points to a critical policy issue: influence rankings alone could potentially prioritise one structural domain over others. From a resource-constrained policy viewpoint, a diversified strategy across clusters may prove more effective than choosing multiple indicators from a single subsystem.

The cluster-based prioritisation strategy therefore shifts the emphasis from ``most influential overall" to ``most influential within structurally distinct pathways." This approach produces a diversified portfolio of entry points, including vaccination coverage (Cluster 1), academic publication output (Cluster 2), and statistical performance (Cluster 3), alongside isolated but synergy-dominated indicators. This diversification improves the ability of interventions to activate multiple independent channels of positive propagation, thereby maximising aggregate system-wide impact while reducing redundancy.

\paragraph{An internally coherent but externally tension-prone subsystem.}

Among the various subsystems, Cluster 4 is particularly notable because of the unique structure it exhibits. This cluster is related to the tightly coupled subsystem with environmental and resource-intensive indicators such as emissions, nitrogen flows, deforestation, pollution burdens, import-related environmental pressures, and governance and social stress indicators. All indicators included in this cluster are trade-off-dominant; however, they are highly positively interlinked with each other. This means that these indicators tend to co-evolve with respect to each other as part of a common structural regime, with changes in any of the environmental and institutional pressures being closely related to changes in the others. The high level of mutual synergy among these indicators means that coordinated policies that focus on improving resource efficiency, reducing environmental pressures such as emissions and pollution, and improving institutional quality could have cumulative effects in this cluster. However, the external trade-off dominance also means that such effects may not automatically translate into socio-economic benefits.

From a policy perspective, this implies that coordinated action within this predominantly environmental and resource-intensive subsystem may yield substantial internal gains, but requires deliberate integration with economic and social policies to mitigate cross-system trade-offs and enhance system-wide coherence.

\paragraph{Implications for national SDG strategy.}

In general, these results indicate that for the successful implementation of the SDG it is essential to focus on interventions that have (i) predominantly positive spillover effects, (ii) significant strength of influence on multiple indicators, and (iii) structurally different propagation pathways.

A purely sectoral strategy may neglect directional spillovers and lead to unforeseen trade-offs. On the other hand, a strategy that focuses on rankings, without accounting for structural clustering, may end up targeting a single subsystem, thereby decreasing the overall diversification of the system. The proposed framework illustrates how a balance-based classification strategy, weighted centrality measures, and flow-based clustering can be combined to produce a more integrated prioritisation strategy.

The structural configuration identified in this analysis suggests several priority areas for policy intervention in India:
\begin{itemize}
    \item Improving primary healthcare and vaccination infrastructure as basic socio-developmental catalysts.
    \item Developing research infrastructure and innovation ecosystems to trigger knowledge-driven spillovers.
    \item Improving statistical and institutional infrastructure to enhance systemic governance feedback.
    \item Framing environmental and infrastructure policies with specific safeguards to avoid trade-offs.
\end{itemize}

\section{CONCLUSIONS}\label{4section:conclusions}
Sustainable development unfolds through interconnected processes rather than isolated sectoral achievements. This research proposes and uses a directed, weighted network approach to analyse how progress spreads through the SDG indicators. Through the analysis of lagged interdependencies, the differentiation of synergies from trade-offs, the measurement of systemic influence, and the detection of structurally coherent paths of propagation, this research goes beyond the assessment of goals to a lag-based and country-specific understanding of development structure.

There are three overarching conclusions that can be drawn.

First, the nature of SDG linkages in the Indian context is highly directional and extends far beyond the boundaries of the goals. Development effects are neither mutually symmetrical nor limited to the institutional setup of the 17-goal framework. Gains in one indicator often precede gains in others, sometimes across sectors, and these effects occur with significant time lags. This reinforces the limitation of silo-based planning and highlights the need for inter-ministerial coordination. Policies formulated without considering directional spillovers are likely to create unintended downstream effects.

Second, although positive spillovers slightly outnumber negative ones, trade-offs remain structurally embedded within the system. Close to one-third of the indicators are trade-off-dominated in their outward spillovers. This result strongly suggests that development outcomes are not universally reinforcing, since gains in areas such as infrastructure, production, and environmental transitions can create adjustment burdens elsewhere. This study does not mean to suggest that such interventions should be avoided. Rather, it draws attention to the need for complementary policies that can absorb or compensate for adverse spillovers. In this regard, trade-offs are not development failures, but indicators of structural tensions that need to be collectively managed.

Third, high-impact interventions are distributed along particular structural paths. Indicators in the categories of health, education, poverty reduction, and governance are shown to have the highest positive downstream impact. However, the rankings of influence alone point to a crucial complexity: the most influential indicators are often spread through intersecting paths. Flow-based clustering illustrates that high levels of synergy are contained within particular subsystems, rather than being evenly distributed throughout the system as a whole. As such, a strategy of choosing top-ranked indicators across multiple subsystems may provide diminishing marginal returns in terms of system-wide diversification.

The methodological contribution of this study lies in integrating complementary elements into a unified framework of prioritisation: (i) the classification of the indicators according to the balance of synergies and trade-offs, (ii) the weighted centrality to estimate the breadth and strength of the influence, and (iii) the flow-based clustering to detect the coherent structures of propagation. These three elements enable a systematic way of identifying development entry points that maximise positive spillovers while taking into account the complexity of the system. It is essential to note that this framework is country-specific and can be applied to other countries with different structural patterns of SDG interactions.

There are a few limitations that must be mentioned. The linear lag-based model framework is used to identify ordered and directional relationships in time and can be seen as a predictive form of causality in the Granger sense. Nevertheless, it does not define structural causal processes or intervention effects. The dominant direction rule adds simplicity by allowing only the dominant directional relationship between two indicators to be considered, which might be a simplification of mutual feedback relationships that could occur at the same time. Future studies could develop this framework using causal inference methods, dynamic simulation models, or cross-country comparisons to analyse how structural SDG patterns change over time.

Overall, the findings suggest that a successful national strategy for achieving the SDGs involves more than simply identifying which goals are important; it also involves understanding how progress moves through the development system. By explicitly modelling direction, strength, and structural propagation, this research shows that sustainable development policy can be made more coherent, less redundant, and more diversified at the structural level. In resource-constrained countries, this type of structurally informed prioritisation may be critical for achieving broad-based progress toward the 2030 Agenda.


\section*{Declaration of Competing Interest}
The authors declare that they have no known competing financial interests or personal relationships that could have appeared to influence the work reported in this paper.

\section*{Funding}
This research was not supported by any external funding.

\section*{Data availability statement}
Data supporting the findings of this study are available from the corresponding author upon reasonable request.

\appendix
\section*{Appendix}\label{Appendix}

\begin{longtable}{clp{11cm}}
\caption{List of SDG Indicators for India}\label{4indicatorslist} \\
\hline
No. & SDG & Indicator Description \\
\hline
\endfirsthead

\hline
No. & SDG & Indicator Description \\
\hline
\endhead

\hline
\endfoot

\hline
\endlastfoot

1 & SDG 1 & Poverty headcount ratio at \$2.15/day (2017 PPP, \%) \\
2 & SDG 1 & Poverty headcount ratio at \$3.65/day (2017 PPP, \%) \\
3 & SDG 2 & Prevalence of undernourishment (\%) \\
4 & SDG 2 & Prevalence of stunting in children under 5 years of age (\%) \\
5 & SDG 2 & Prevalence of wasting in children under 5 years of age (\%) \\
6 & SDG 2 & Prevalence of obesity, BMI $\geq 30$ (\% of adult population) \\
7 & SDG 2 & Human Trophic Level (best 2-3 worst) \\
8 & SDG 2 & Cereal yield (tonnes per hectare of harvested land) \\
9 & SDG 2 & Sustainable Nitrogen Management Index (best 0-1.41 worst)  \\
10 & SDG 3 & Maternal mortality rate (per 100,000 live births) \\
11 & SDG 3 & Neonatal mortality rate (per 1,000 live births) \\
12 & SDG 3 & Mortality rate, under-5 (per 1,000 live births) \\
13 & SDG 3 & Incidence of tuberculosis (per 100,000 population) \\
14 & SDG 3 & Age-standardized death rate due to cardiovascular disease, cancer, diabetes, or chronic respiratory disease in adults aged 30–70 years (\%) \\
15 & SDG 3 & Traffic deaths (per 100,000 population) \\
16 & SDG 3 & Life expectancy at birth (years) \\
17 & SDG 3 & Adolescent fertility rate (births per 1,000 females aged 15 to 19)  \\
18 & SDG 3 & Births attended by skilled health personnel (\%) \\
19 & SDG 3 & Surviving infants who received 2 WHO-recommended vaccines (\%) \\
20 & SDG 3 & Universal health coverage (UHC) index of service coverage (worst 0-100 best) \\
21 & SDG 3 & Subjective well-being (average ladder score, worst 0-10 best) \\
22 & SDG 4 & Participation rate in pre-primary organized learning (\% of children aged 4 to 6) \\
23 & SDG 4 & Net primary enrollment rate (\%) \\
24 & SDG 4 & Lower secondary completion rate (\%) \\
25 & SDG 4 & Literacy rate (\% of population aged 15 to 24) \\
26 & SDG 5 & Demand for family planning satisfied by modern methods (\% of females aged 15 to 49) \\
27 & SDG 5 & Ratio of female-to-male mean years of education received (\%) \\
28 & SDG 5 & Ratio of female-to-male labor force participation rate (\%) \\
29 & SDG 5 & Seats held by women in national parliament (\%) \\
30 & SDG 6 & Population using at least basic drinking water services (\%) \\
31 & SDG 6 & Population using at least basic sanitation services (\%) \\
32 & SDG 6 & Freshwater withdrawal (\% of available freshwater resources) \\
33 & SDG 6 & Scarce water consumption embodied in imports ($m^3H_2O$ eq/capita) \\
34 & SDG 7 & Population with access to electricity (\%) \\
35 & SDG 7 & Population with access to clean fuels and technology for cooking (\%) \\
36 & SDG 7 & $CO_2$ emissions from fuel combustion per total electricity output $(MtCO_2/TWh)$  \\
37 & SDG 7 & Renewable energy share in total final energy consumption (\%) \\
38 & SDG 8 & Adults with an account at a bank or other financial institution or with a mobile-money-service provider (\% of population aged 15 or over) \\
39 & SDG 8 & Unemployment rate (\% of total labor force, ages 15+) \\
40 & SDG 8 & Fundamental labor rights are effectively guaranteed (worst 0–1 best) \\
41 & SDG 8 & Fatal work-related accidents embodied in imports (per million population) \\
42 & SDG 9 & Population using the internet (\%) \\
43 & SDG 9 & Mobile broadband subscriptions (per 100 population) \\
44 & SDG 9 & Logistics Performance Index: Infrastructure Score (worst 1–5 best) \\
45 & SDG 9 & The Times Higher Education Universities Ranking: Average score of top 3 universities (worst 0-100 best) \\
46 & SDG 9 & Articles published in academic journals (per 1,000 population) \\
47 & SDG 9 & Expenditure on research and development (\% of GDP) \\
48 & SDG 9 & Total patent applications by applicant's origin (per million population) \\
49 & SDG 10 & Gini coefficient \\
50 & SDG 10 & Palma ratio \\
51 & SDG 11 & Proportion of urban population living in slums (\%) \\
52 & SDG 11 & Annual mean concentration of PM2.5 ($\mu g / m^3$) \\
53 & SDG 11 & Access to improved water source, piped (\% of urban population) \\
54 & SDG 12 & Production-based air pollution (DALYs per 1,000 population) \\
55 & SDG 12 & Air pollution associated with imports (DALYs per 1,000 population) \\
56 & SDG 12 & Production-based nitrogen emissions (kg/capita) \\
57 & SDG 12 & Nitrogen emissions associated with imports (kg/capita) \\
58 & SDG 13 & $CO_2$ emissions from fossil fuel combustion and cement production (tCO2/capita) \\
59 & SDG 13 & GHG emissions embodied in imports ($tCO_2$/capita) \\
60 & SDG 14 & Ocean Health Index: Clean Waters score (worst 0-100 best) \\
61 & SDG 14 & Fish caught from overexploited or collapsed stocks (\% of total catch) \\
62 & SDG 14 & Fish caught by trawling or dredging (\%) \\
63 & SDG 14 & Fish caught that are then discarded (\%) \\
64 & SDG 15 & Mean area that is protected in terrestrial sites important to biodiversity (\%) \\
65 & SDG 15 & Mean area that is protected in freshwater sites important to biodiversity (\%) \\
66 & SDG 15 & Red List Index of species survival (worst 0-1 best)  \\
67 & SDG 15 & Permanent deforestation (\% of forest area, 3-year average) \\
68 & SDG 15 & Imported deforestation ($m^2$/capita) \\
69 & SDG 16 & Homicides (per 100,000 population) \\
70 & SDG 16 & Crime is effectively controlled \\
71 & SDG 16 & Unsentenced detainees (\% of prison population) \\
72 & SDG 16 & Corruption Perceptions Index (worst 0-100 best) \\
73 & SDG 16 & Press Freedom Index (worst 0-100 best) \\
74 & SDG 16 & Access to and affordability of justice (worst 0–1 best) \\
75 & SDG 16 & Timeliness of administrative proceedings (worst 0 - 1 best) \\
76 & SDG 16 & Expropriations are lawful and adequately compensated (worst 0 - 1 best) \\
77 & SDG 17 & Government spending on health and education (\% of GDP) \\
78 & SDG 17 & Other countries: Government revenue excluding grants (\% of GDP) \\
79 & SDG 17 & Statistical Performance Index (worst 0-100 best) \\

\end{longtable}

\bibliographystyle{apalike}
\bibliography{references}

@article{kottari2026network,
  title={A Network-Based Framework to Identify Synergies and Trade offs among SDG Indicators},
  author={Kottari, Gaurav and Azhad, Qazi J and Sahni, Niteesh},
  journal={arXiv preprint arXiv:2602.01676},
  year={2026}
}

@article{albrizio2014empirical,
  title={Empirical evidence on the effects of environmental policy stringency on productivity growth},
  author={Albrizio, Silvia and Kozluk, Tomasz and Zipperer, Vera},
  journal={OECD Economic Department Working Papers},
  number={1179},
  pages={0\_1},
  year={2014},
  publisher={Organisation for Economic Cooperation and Development (OECD)}
}

@article{cao2017effect,
  title={The effect of environmental regulation on employment in resource-based areas of China—an empirical research based on the mediating effect model},
  author={Cao, Wenbin and Wang, Hui and Ying, Huihui},
  journal={International journal of environmental research and public health},
  volume={14},
  number={12},
  pages={1598},
  year={2017},
  publisher={MDPI}
}

@article{zhao2021quantifying,
  title={Quantifying economic-social-environmental trade-offs and synergies of water-supply constraints: An application to the capital region of China},
  author={Zhao, Dandan and Liu, Junguo and Sun, Laixiang and Ye, Bin and Hubacek, Klaus and Feng, Kuishuang and Varis, Olli},
  journal={Water research},
  volume={195},
  pages={116986},
  year={2021},
  publisher={Elsevier}
}

@article{bloom2004effect,
  title={The effect of health on economic growth: a production function approach},
  author={Bloom, David E and Canning, David and Sevilla, Jaypee},
  journal={World development},
  volume={32},
  number={1},
  pages={1--13},
  year={2004},
  publisher={Elsevier}
}

@article{pradhan2017systematic,
  title={A systematic study of sustainable development goal (SDG) interactions},
  author={Pradhan, Prajal and Costa, Lu{\'\i}s and Rybski, Diego and Lucht, Wolfgang and Kropp, J{\"u}rgen P},
  journal={Earth's future},
  volume={5},
  number={11},
  pages={1169--1179},
  year={2017},
  publisher={Wiley Online Library}
}

@article{song2023unpacking,
  title={Unpacking the sustainable development goals (SDGs) interlinkages: A semantic network analysis of the SDGs targets},
  author={Song, Jaemin and Jang, Chang-Ho},
  journal={Sustainable development},
  volume={31},
  number={4},
  pages={2784--2796},
  year={2023},
  publisher={Wiley Online Library}
}

@article{swain2021modeling,
  title={Modeling interlinkages between sustainable development goals using network analysis},
  author={Swain, Ranjula Bali and Ranganathan, Shyam},
  journal={World Development},
  volume={138},
  pages={105136},
  year={2021},
  publisher={Elsevier}
}

@article{nilsson2016policy,
  title={Policy: map the interactions between sustainable development goals},
  author={Nilsson, M{\aa}ns and Griggs, Dave and Visbeck, Martin},
  journal={Nature},
  volume={534},
  number={7607},
  pages={320--322},
  year={2016},
  publisher={Nature Publishing Group UK London}
}

@article{le2015towards,
  title={Towards integration at last? The sustainable development goals as a network of targets},
  author={Le Blanc, David},
  journal={Sustainable Development},
  volume={23},
  number={3},
  pages={176--187},
  year={2015},
  publisher={Wiley Online Library}
}

@article{assembly2015resolution,
  title={Resolution adopted by the General Assembly on 11 September 2015},
  author={Assembly, General},
  journal={New York: United Nations},
  pages={14},
  year={2015}
}

@article{allen2019prioritising,
  title={Prioritising SDG targets: Assessing baselines, gaps and interlinkages},
  author={Allen, Cameron and Metternicht, Graciela and Wiedmann, Thomas},
  journal={Sustainability Science},
  volume={14},
  pages={421--438},
  year={2019},
  publisher={Springer}
}

@article{weitz2018towards,
  title={Towards systemic and contextual priority setting for implementing the 2030 Agenda},
  author={Weitz, Nina and Carlsen, Henrik and Nilsson, M{\aa}ns and Sk{\aa}nberg, Kristian},
  journal={Sustainability science},
  volume={13},
  pages={531--548},
  year={2018},
  publisher={Springer}
}

@article{khot2026existing,
  title={Existing gaps in understanding Sustainable Development Goals interactions: Insights from a systematic review},
  author={Khot, Utkarsh Ashok and Warchold, Anne and Pradhan, Prajal},
  journal={Environmental Impact Assessment Review},
  volume={118},
  pages={108274},
  year={2026},
  publisher={Elsevier}
}

@article{opsahl2010node,
  title={Node centrality in weighted networks: Generalizing degree and shortest paths},
  author={Opsahl, Tore and Agneessens, Filip and Skvoretz, John},
  journal={Social networks},
  volume={32},
  number={3},
  pages={245--251},
  year={2010},
  publisher={Elsevier}
}

@book{griggs2017guide,
  title={A guide to SDG interactions: from science to implementation},
  author={Griggs, DJ and Nilsson, Mans and Stevance, A and McCollum, David and others},
  year={2017},
  publisher={International Council for Science, Paris}
}

@article{pham2020interactions,
  title={Interactions among Sustainable Development Goals: Knowledge for identifying multipliers and virtuous cycles},
  author={Pham-Truffert, Myriam and Metz, Florence and Fischer, Manuel and Rueff, Henri and Messerli, Peter},
  journal={Sustainable development},
  volume={28},
  number={5},
  pages={1236--1250},
  year={2020},
  publisher={Wiley Online Library}
}

@article{tremblay2020sustainable,
  title={Sustainable development goal interactions: An analysis based on the five pillars of the 2030 agenda},
  author={Tremblay, David and Fortier, Fran{\c{c}}ois and Boucher, Jean-Fran{\c{c}}ois and Riffon, Olivier and Villeneuve, Claude},
  journal={Sustainable Development},
  volume={28},
  number={6},
  pages={1584--1596},
  year={2020},
  publisher={Wiley Online Library}
}

@article{wong2019avoidance,
  title={Avoidance of conflicts and trade-offs: A challenge for the policy integration of the United Nations Sustainable Development Goals},
  author={Wong, Ryan and van der Heijden, Jeroen},
  journal={Sustainable Development},
  volume={27},
  number={5},
  pages={838--845},
  year={2019},
  publisher={Wiley Online Library}
}

@article{xiao2023synergies,
  title={Synergies and trade-offs across sustainable development goals: A novel method incorporating indirect interactions analysis},
  author={Xiao, Huijuan and Liu, Yue and Ren, Jingzheng},
  journal={Sustainable Development},
  volume={31},
  number={2},
  pages={1135--1148},
  year={2023},
  publisher={Wiley Online Library}
}

@article{fader2018toward,
  title={Toward an understanding of synergies and trade-offs between water, energy, and food SDG targets},
  author={Fader, Marianela and Cranmer, Colleen and Lawford, Richard and Engel-Cox, Jill},
  journal={Frontiers in Environmental Science},
  volume={6},
  pages={112},
  year={2018},
  publisher={Frontiers Media SA}
}

@article{de2020synergies,
  title={Synergies and trade-offs among sustainable development goals: the case of Spain},
  author={de Miguel Ramos, Carlos and Laurenti, Rafael},
  journal={Sustainability},
  volume={12},
  number={24},
  pages={10506},
  year={2020},
  publisher={MDPI}
}

@article{kostetckaia2022sustainable,
  title={How sustainable development goals interlinkages influence European Union countries’ progress towards the 2030 agenda},
  author={Kostetckaia, Mariia and Hametner, Markus},
  journal={Sustainable Development},
  volume={30},
  number={5},
  pages={916--926},
  year={2022},
  publisher={Wiley Online Library}
}

@misc{hegre2020synergies,
  title={Synergies and Trade-Offs in Reaching the Sustainable Development Goals. Sustainability, 12 (20), 8729},
  author={Hegre, H and Petrova, K and von Uexkull, N},
  year={2020}
}

@article{miao2025priority,
  title={The Priority Implementation of the Sustainable Development Goals in Underdeveloped Mountain Regions},
  author={Miao, Junxia and Song, Xiaoyu and Zhong, Fanglei and Gao, Feng and Huang, Chunlin and Zhao, Xueyan},
  journal={Sustainable Development},
  volume={33},
  number={3},
  pages={3331--3347},
  year={2025},
  publisher={Wiley Online Library}
}

@article{griggs2013sustainable,
  title={Sustainable development goals for people and planet},
  author={Griggs, David and Stafford-Smith, Mark and Gaffney, Owen and Rockstr{\"o}m, Johan and {\"O}hman, Marcus C and Shyamsundar, Priya and Steffen, Will and Glaser, Gisbert and Kanie, Norichika and Noble, Ian},
  journal={Nature},
  volume={495},
  number={7441},
  pages={305--307},
  year={2013},
  publisher={Nature Publishing Group UK London}
}

@article{stafford2017integration,
  title={Integration: the key to implementing the Sustainable Development Goals},
  author={Stafford-Smith, Mark and Griggs, David and Gaffney, Owen and Ullah, Farooq and Reyers, Belinda and Kanie, Norichika and Stigson, Bjorn and Shrivastava, Paul and Leach, Melissa and O’Connell, Deborah},
  journal={Sustainability science},
  volume={12},
  number={6},
  pages={911--919},
  year={2017},
  publisher={Springer}
}

@article{biermann2017global,
  title={Global governance by goal-setting: the novel approach of the UN Sustainable Development Goals},
  author={Biermann, Frank and Kanie, Norichika and Kim, Rakhyun E},
  journal={Current opinion in environmental sustainability},
  volume={26},
  pages={26--31},
  year={2017},
  publisher={Elsevier}
}

@article{tosun2017governing,
  title={Governing the interlinkages between the sustainable development goals: Approaches to attain policy integration},
  author={Tosun, Jale and Leininger, Julia},
  journal={Global challenges},
  volume={1},
  number={9},
  pages={1700036},
  year={2017},
  publisher={Wiley Online Library}
}

@article{cejudo2017addressing,
  title={Addressing fragmented government action: Coordination, coherence, and integration},
  author={Cejudo, Guillermo M and Michel, Cynthia L},
  journal={Policy Sciences},
  volume={50},
  number={4},
  pages={745--767},
  year={2017},
  publisher={Springer}
}

@article{alcamo2020analysing,
  title={Analysing interactions among the sustainable development goals: findings and emerging issues from local and global studies},
  author={Alcamo, Joseph and Thompson, John and Alexander, Anthony and Antoniades, Andreas and Delabre, Izabela and Dolley, Jonathan and Marshall, Fiona and Menton, Mary and Middleton, Jo and Scharlemann, J{\"o}rn PW},
  journal={Sustainability Science},
  volume={15},
  number={6},
  pages={1561--1572},
  year={2020},
  publisher={Springer}
}

@article{kroll2019sustainable,
  title={Sustainable Development Goals (SDGs): Are we successful in turning trade-offs into synergies?},
  author={Kroll, Christian and Warchold, Anne and Pradhan, Prajal},
  journal={Palgrave Communications},
  volume={5},
  number={1},
  year={2019},
  publisher={Springer Science and Business Media LLC}
}

@article{sachs2019six,
  title={Six transformations to achieve the sustainable development goals},
  author={Sachs, Jeffrey D and Schmidt-Traub, Guido and Mazzucato, Mariana and Messner, Dirk and Nakicenovic, Nebojsa and Rockstr{\"o}m, Johan},
  journal={Nature sustainability},
  volume={2},
  number={9},
  pages={805--814},
  year={2019},
  publisher={Nature Publishing Group UK London}
}

@article{breuer2019translating,
  title={Translating sustainable development goal (SDG) interdependencies into policy advice},
  author={Breuer, Anita and Janetschek, Hannah and Malerba, Daniele},
  journal={Sustainability},
  volume={11},
  number={7},
  pages={2092},
  year={2019},
  publisher={MDPI}
}

@article{bennich2020deciphering,
  title={Deciphering the scientific literature on SDG interactions: A review and reading guide},
  author={Bennich, Therese and Weitz, Nina and Carlsen, Henrik},
  journal={Science of the Total Environment},
  volume={728},
  pages={138405},
  year={2020},
  publisher={Elsevier}
}

@article{scharlemann2020towards,
  title={Towards understanding interactions between Sustainable Development Goals: The role of environment--human linkages},
  author={Scharlemann, J{\"o}rn PW and Brock, Rebecca C and Balfour, Nicholas and Brown, Claire and Burgess, Neil D and Guth, Miriam K and Ingram, Daniel J and Lane, Richard and Martin, Juliette GC and Wicander, Sylvia and others},
  journal={Sustainability science},
  volume={15},
  number={6},
  pages={1573--1584},
  year={2020},
  publisher={Springer}
}

@article{rosvall2008maps,
  title={Maps of random walks on complex networks reveal community structure},
  author={Rosvall, Martin and Bergstrom, Carl T},
  journal={Proceedings of the national academy of sciences},
  volume={105},
  number={4},
  pages={1118--1123},
  year={2008},
  publisher={National Academy of Sciences}
}

@article{kullmann2002time,
  title={Time-dependent cross-correlations between different stock returns: A directed network of influence},
  author={Kullmann, L{\'a}szl{\'o} and Kert{\'e}sz, Janos and Kaski, Kimmo},
  journal={Physical Review E},
  volume={66},
  number={2},
  pages={026125},
  year={2002},
  publisher={APS}
}

@article{li2022undirected,
  title={Undirected and directed network analysis of the chinese stock market},
  author={Li, Binghui and Yang, Yuehan},
  journal={Computational Economics},
  volume={60},
  number={3},
  pages={1155--1173},
  year={2022},
  publisher={Springer}
}

@article{collste2017policy,
  title={Policy coherence to achieve the SDGs: using integrated simulation models to assess effective policies},
  author={Collste, David and Pedercini, Matteo and Cornell, Sarah E},
  journal={Sustainability science},
  volume={12},
  number={6},
  pages={921--931},
  year={2017},
  publisher={Springer}
}

@misc{mapequation2026software,
    title = {{The MapEquation software package}},
    author = {Edler, Daniel and Holmgren, Anton and Rosvall, Martin},
    howpublished = {\url{https://mapequation.org}},
    year = 2026,
}

@book{sachs2025sdr,
  author    = {Sachs, Jeffrey D. and Lafortune, Guillaume and Fuller, Grayson and Iablonovski, Guilherme},
  title     = {Financing Sustainable Development to 2030 and Mid-Century. Sustainable Development Report 2025},
  year      = {2025},
  publisher = {Sustainable Development Solutions Network (SDSN) and Dublin University Press},
  address   = {Paris and Dublin},
  doi       = {10.25546/111909},
  url       = {https://doi.org/10.25546/111909}
}

\end{document}